\documentclass[a4paper,12pt]{amsart}
\usepackage{amssymb,amsmath,amsthm,amscd}
\newcommand{\rn}{{\mathbf R}^n}
\newcommand{\ia}{I_\alpha}

\begin{document}
\title{Morrey spaces and fractional integral operators}

\author[Eridani, Kokilashvili,]{Eridani$^\diamond$, Vakhtang Kokilashvili,$^\dagger$}
\author[Meskhi]{Alexander Meskhi$^\ddagger$\\\\}


\maketitle

\vspace{-1.5cm}

{\begin{center}
{\small\itshape{$^\diamond$Department of Mathematics, \\
Airlangga University, Campus C, Mulyorejo,\\
Surabaya 60115, Indonesia. \\
\email{E-mail address: keumala\_ikhsanti@yahoo.com.}\\

\bigskip

$^\dagger$$^{,\ddagger}$A. Razmadze Mathematical Institute,\\
M. Aleksidze St. 0193, Tbilisi 380093, Georgia. \\
\email{E-mail address:
\{$^\ddagger$meskhi,$^\dagger$kokil\,\}@rmi.acnet.ge.}}}

\bigskip

$^\diamond$$^{,\ddagger}${\it School of Mathematical Sciences,\\
Government College University,\\
68-B New Muslim Town, Lahore, Pakistan.} \\
\end{center}}

\begin{abstract}
The present paper is devoted to the boundedness of fractional integral operators in Morrey spaces defined on quasimetric measure spaces. In particular, Sobolev, trace and weighted inequalities with power weights for potential operators are established. In the case when measure satisfies the doubling condition the derived conditions are simultaneously necessary and sufficient for appropriate inequalities.
\end{abstract}

\vskip+0.5cm


\vskip+0.5cm

\vskip+0.5cm

\section{Introduction}
The main purpose of this paper is to establish the boundedness  of
fractional integral operators in (weighted) Morrey spaces defined
on quasimetric measure spaces. We derive Sobolev, trace and two-weight inequalities for fractional
integrals. In particular, we generalize: a) D. Adams \cite{Adm} trace inequality;
b) the theorem by E. M. Stein and G. Weiss  \cite{ST1} regarding the two-weight inequality for the Riesz potentials; c) Sobolev-type
inequality. We emphasize that in the most cases the derived conditions are necessary and sufficient for appropriate inequalities.
\medskip

In the paper \cite{GGK} (see also \cite{GGKK}, Ch. 2) integral-type sufficient condition
guaranteeing the two-weight weak-type inequality for integral
operator with positive kernel defined on nonhomogeneous spaces was
established. In the same paper (see also \cite{GGKK}, Ch. 2) the authors solved the two-weight problem for
kernel operators on spaces of homogeneous type.

In \cite{KM0} (see also \cite{EKM}, Ch.6) a complete description
of non-doubling measure $\mu$ guaranteeing the boundedness of
fractional integral operator $\ia$ (see the next section for the
definition) from $L^p(\mu, X)$ to $L^q(\mu, X),\quad
1<p<q<\infty,$ was given. We notice that this result for potentials was derived in \cite{Ko} for potentials on Euclidean spaces. In \cite{KM0}, theorems of Sobolev and Adams type for fractional integrals defined on
quasimetric measure spaces were established. For the boundedness of fractional integrals on metric measure spaces we refer also to \cite{GCG}. Some
two-weight norm inequalities for fractional operators on $\rn$ with non-doubling measure were studied in \cite{GCM}. Further, in the paper \cite{KM1} necessary and sufficient conditions on measure $\mu$ governing the inequality of Stein-Weiss type on nonhomogeneous spaces were established.

The boundedness of the Riesz potential in Morrey spaces defined on
Euclidean spaces was studied in \cite{Pe} and \cite{Ad1}. The same
problem for fractional integrals on ${\mathbb{R}}^n$ with
non-doubling measure was investigated in \cite{ST}.

Finally we mention that necessary and sufficient conditions for the boundedness of maximal
operators and Riesz potentials in the local Morrey-type spaces were derived in
\cite{BG}, \cite{BGG}.

The main results of this paper were presented in \cite{ErKoMe}.

It should be emphasized that the results of this work are new even
for Euclidean spaces.

Constants (often different constants in the same series of inequalities)
will generally be denoted by $c$ or $C$.

\section{Preliminaries}
Throughout the paper we assume that $X:=(X, \rho, \mu)$ is a
topological space, endowed with a complete measure $\mu$ such that
the space of compactly supported continuous functions is dense in
$L^1(X,\mu)$ and there exists a function (quasimetric) $\rho
:X\times X \longrightarrow [0, \infty)$ satisfying the conditions:

\medskip

\noindent(1) $\rho(x,y)>0$ for all $x\neq y,$ and $\rho(x,x)=0$ for all $x\in X$;

\noindent(2) there exists a constant $a_0\geq 1$, such that
$\rho(x,y)\leq a_0\rho(y,x)$ for all $x,\,y\in X;$

\noindent(3) there exists a constant $a_1\geq 1$, such that
$\rho(x,y)\leq a_1(\rho(x,z)+\rho(z,y))$ for all $x,\,y,\,z\in X.$

\medskip

\noindent We assume that the balls $B(a,r):=\{x\in
X:\rho(x,a)<r\}$ are measurable, for $a \in X,\,r>0,$ and
$0\leq\mu (B(a,r))<\infty.$ For every neighborhood $V$ of $x\in
X,$ there exists $r>0,$ such that $B(x,r)\subset V.$ We also
assume that $\mu(X)=\infty,\,\mu\{a\}=0,$ and $B(a,r_2)\setminus
B(a,r_1)\neq \emptyset,$ for all $a\in X,\,\,0<r_1<r_2<\infty.$

The triple $(X, \rho, \mu)$ will be called quasimetric measure space.
\medskip

Let $0<\alpha<1$. We consider the fractional integral operators $\ia,$ and
$K_{\alpha}$ given by
$$\ia f(x):=\int_X
f(y){\rho(x,y)}^{\alpha-1}\,d\mu(y), $$
$$ K_{\alpha}f(x):=\int_X f(y){(\mu B(x,\rho(x,y)))}^{\alpha-1}\,d\mu(y), $$
for suitable $f$ on $X$.
\medskip

Suppose that $\nu$ is another measure on $X$,  $\lambda\geq 0$ and $1\leq p<\infty$. We deal with the Morrey space
$L^{p,\lambda}(X, \nu,\mu)$, which is the set of all functions $f\in L^p_{\rm loc}(X,\nu)$ such that
$$ {\|f\|}_{L^{p,\lambda}(X,\nu,\mu)}:=\sup_{B}{\left( \frac{1}{\mu (B)^{\lambda}} \int_{B}{|f(y)|}^p \,d\nu(y)\right)
}^{1/p}<\infty,  $$
where the supremum is taken over all balls $B$.

If $\nu=\mu$, then we have the classical
Morrey space $L^{p,\lambda}(X,\mu)$ with measure $\mu$. When
$\nu=\mu$ and $\lambda =0$, then $L^{p,\lambda}(X,\nu,\mu)
=L^p(X,\mu)$ is the Lebesgue space with measure $\mu$.

Further, suppose that $\beta \in \mathbf{R}.$ We are also
interested in weighted Morrey space $M_\beta^{p,\lambda}(X,\mu)$
which is the set of all $\mu$-measurable functions $f$ such that
$$ {\|f\|}_{M_\beta^{p,\lambda}(X,\mu)}:=\sup\limits_{a\in X; r>0}{\left( \frac{1}{r^{\lambda}}\int_{B(a,r)}{|f(y)|}^p\rho(a,y)^\beta
\,d\mu(y)\right) }^{1/p}<\infty. $$

If $\beta=0$, then we denote $M^{p,\lambda}_{\beta}(X,\mu):= M^{p,\lambda}(X,\mu)$.

\medskip

We say that a measure $\mu$  satisfies  the growth condition ($\mu
\in (GC)$), if there exists $C_0>0$ such that $\mu(B(a,r))\leq C_0
r$; further,  $\mu$ satisfies the doubling condition ($\mu \in$
(DC)) if $\mu(B(a,2r))\leq C_1\,\mu(B(a,r))$ for some $C_1>1.$ If
$\mu \in (DC)$, then $(X, \rho, \mu)$ is called a space of
homogeneous type (SHT). A quasimetric measure space $(X, \rho,
\mu)$, where the doubling condition might be failed, is also
called a non-homogeneous space.

The measure $\mu$ on $X$ satisfies the reverse
doubling condition ($\mu \in (RDC)$) if there are  constants
$\eta_1$ and $\eta_2$ with $\eta_1>1$ and $\eta_2>1$ such that
$$ \mu B(x, \eta_1 r) \geq \eta_2 \mu B(x, r). \eqno{(1)}$$


\vskip+0.1cm

It is known (see e.g. \cite{StTo}, p. 11) that if $\mu \in \; (DC)$, then $\mu \in \; (RDC)$.

\medskip

The next statements is from \cite{KM0} (see also \cite[Theorem
6.1.1, Corollary 6.1.1]{EKM} and \cite{Ko} in the case of
Euclidean spaces).

\bigskip

{\bf Theorem A.} {\em Let $(X, \rho, \mu)$ be a quasimetric measure space. Suppose that $1<p<q<\infty$ and $0<\alpha<1$. Then $\ia$ is bounded
from $L^p(X)$ to $L^q(X)$ if and only if there exists a positive constant $C$  such
that $$  \mu(B(a,r))\le C r^{s},\quad
s=\frac{pq(1-\alpha)}{pq+p-q}, \eqno{(2)}$$
for all $a\in X$ and $r>0$.}

\vskip+0.1cm

{\bf Corollary B.} {\em  Let  Let $(X, \rho, \mu)$ be a quasimetric measure space, $1<p<1/\alpha$ and $1/q=1/p-\alpha.$
Then $\ia$ is bounded from $L^p(X)$ to $L^q(X)$ if and only if $\mu \in$ (GC). }

\bigskip

The latter statement by different proof was also derived in
\cite{GCG} for metric spaces.

We to prove some of our statements we need the following
Hardy-type transform:

$$ H_a f(x):=\int_{\rho(a,y)\le \rho(a,x)}f(y)\,d\mu(y), $$
where $a$ is a fixed point of $X$ and $f\in L_{\rm loc}(X,\mu)$.

\medskip


{\bf Theorem C.} {\em Suppose that $(X, \rho, \mu)$ be a quasimetric measure space, $1<p\le q<\infty$ and  $V$ and $W$
are non-negative functions defined on $X\times X$. Let $\nu$ be another measures on $X$. If there exists  a positive constant $C$ such that for every
$a\in X$ and $t>0$,
$$
\left(\int_{\rho(a,y)\geq t}V(a,y)\,d\nu(y)\right)^{1/q}
\left(\int_{\rho(a,y)\le t}W(a,y)^{1-p'}\,d\mu(y)\right)^{1/p'}
\le C<\infty, $$ then there exists a positive constant $c$ such
that for all $f\geq 0$ and $a\in X$ the inequality
$$ \left(\int_{B(a,r)} (H_af(x))^qV(a,x)\,d\nu(x)\right)^{1/q}\le
c\,\left(\int_{B(a,r)}(f(x))^pW(a,x)\,d\mu(x)\right)^{1/p} $$
holds. }

\bigskip

This statement was proved in \cite[Section 1.1]{EKM} for Lebesgue
spaces.

\bigskip

{\em  Proof of Theorem C.} Let $f\geq 0$. We define
$S(s):=\int_{\rho(a,y)<s} f(y)\,d\mu(y),$ for $s\in [0,r].$
Suppose $S(r)<\infty,$ then $2^m<S(r)\le 2^{m+1},$ for some $m\in
\mathbb{Z}.$ Let
$$ s_j:=\sup \{t: S(t)\le 2^j\},\,\,j\le m,\quad {\rm and}\quad
s_{m+1}:=r. $$
Then it is easy to see that (see also \cite[pp.\,5-8]{EKM} for details)
$(s_j)^{m+1}_{j=-\infty}$ is a non-decreasing sequence, $S(s_j)\le
2^j,\,S(t)\geq 2^j$ for $t>s_j,$ and
$$
2^j\le \int_{s_j\le \rho(a,y)\le s_{j+1}} f(y)\,d\mu(y).
$$
If $\beta:=\lim\limits_{j\rightarrow -\infty}s_j,$ then
$$
\rho(a,x)<r \Leftrightarrow \rho(a,x)\in [0,\beta]\cup
\bigcup_{j=-\infty}^m(s_j,s_{j+1}].
$$
If $S(r)=\infty,$ then we may put $m=\infty.$ Since
$$
0\le \int_{\rho(a,y)<\beta}f(y)\,d\mu(y)\le S(s_j)\le 2^j,
$$
for every $j,$ therefore  $\int_{\rho(a,y)<\beta}f(y)\,d\mu(y)=0.$ From
these observations, we have

$$ \int_{\rho(a,x)<r}(H_af(x))^q V(a,x)\,d\nu(x) \le
\sum_{j=-\infty}^m \int_{s_j\le \rho(a,x)\le s_{j+1}}(H_af(x))^q
V(a,x)\,d\nu(x) $$

$$ \le \sum_{j=-\infty}^m \int_{s_j\le \rho(a,x)\le s_{j+1}}V(a,x)
\left(\int_{\rho(a,y)\le s_{j+1}}(f(y))\,d\mu(y)\right)^qd\nu(x).
$$

Notice that

$$ \int_{\rho(a,y)\le s_{j+1}} f d\mu \leq S(s_{j+2}) \le 2^{j+2}\le C\,\int_{s_{j-1}\le \rho(a,y)\le
s_j}fd\mu. $$

Using H\"{o}lder's inequality, we find that

$$ \int_{\rho(a,x)<r}(H_af(x))^qV(a,x)\,d\mu(x) $$

$$ \le \sum_{j=-\infty}^m \int_{s_j\le \rho(a,x)\le s_{j+1}} V(a,x)
\left(\int_{\rho(a,y)\le s_{j+1}}(f(y))\,d\mu(y)\right)^qd\nu(x)$$

$$ \le C\,\sum_{j=-\infty}^m \int_{s_j\le \rho(a,x)\le s_{j+1}}
V(a,x) \left(\int_{s_{j-1}\leq \rho(a,y)\le
s_j}(f(y))\,d\mu(y)\right)^q d\nu(x) $$

$$ \le C\,\sum_{j=-\infty}^m \int_{s_j\le \rho(a,x)\le
s_{j+1}}V(a,x)\,d\nu(x)\left(\int_{s_{j-1}\leq \rho(a,y)\le
s_j}(f(y))^pW(a,y)\,d\mu(y)\right)^{q/p} $$

$$ \times \left(\int_{s_{j-1}\leq \rho(a,y)\le
s_j}W(a,y)^{1-p'}\,d\mu(y)\right)^{q/p'} $$

$$ \le C\,\sum_{j=-\infty}^m \int_{s_j\le
\rho(a,y)}V(a,y)\,d\nu(y)\left(\int_{\rho(a,y)\le
s_j}W(a,y)^{1-p'}\,d\mu(y)\right)^{q/p'} $$

$$ \times \left(\int_{s_{j-1}\leq \rho(a,y)\le
s_j}(f(y))^p W(a,y)\,d\mu(y)\right)^{q/p}$$

$$ \le C\, \sum_{j=-\infty}^m \left(\int_{s_{j-1}\leq \rho(a,y)\le
s_j}(f(y))^pW(a,y)\,d\mu(y)\right)^{q/p}$$

$$ \le C\,\left(\int_{\rho(a,y)\le
r}(f(y))^pW(a,y)\,d\mu(y)\right)^{q/p}. $$

This completes the proof of the theorem. $\;\;\;\; \Box$

\bigskip

For our purposes we also need the following lemma (see \cite{KM2}
for the case of $\rn$).

\bigskip
{\bf Lemma D.} {\em Suppose that $(X, \rho, \mu)$ be an $SHT$. Let $0<\lambda<1\leq p<\infty$.  Then there exists
a positive constant $C$ such that for all balls $B_0$,

$$ {\|\chi_{B_0}\|}_{L^{p,\lambda}(X,\mu)}\le C {\mu(B_0)}^{(1-\lambda)/p}. $$ }

\medskip

{\em Proof.} Let $B_0:=B(x_0,r_0)$ and $B:=B(a,r).$ We have

$$
{\|\chi_{B_0}\|}_{L^{p,\lambda}(X, \mu)} = \sup\limits_{B}{\left(
\frac{\mu(B_0\cap B)}{{\mu(B)}^{\lambda}}\right) }^{1/p}. $$

Suppose that $B_0\cap B\neq \emptyset$. Let us assume that $r\leq r_0.$ Then (see \cite{StTo}, Lemma 1, or \cite{GGKK}, p.9)
$ B\subset B(x_0, br_0)$, where $b= a_1(1+a_0)$. By the doubling condition it follows that
$$ \frac{\mu(B\cap B_0)}{\mu(B)^{\lambda}}
\leq \frac{\mu(B)}{\mu(B)^\lambda} =\mu(B)^{1-\lambda}  \leq
\mu(B(x_0, br_0))^{1-\lambda} $$
$$\leq C\,\mu(B_0)^{1-\lambda}.$$

Let now $r_0<r$. Then $\mu B_0\leq c \mu B$, where the constant
$c$ depends only on $a_1$ and $a_0$. Then

$$ \frac{\mu(B \cap B_0)}{\mu(B)^{\lambda}}
\leq c \frac{\mu(B_0)}{\mu(B_0)^\lambda} = c \mu
(B_0)^{1-\lambda}. $$

$\Box$

\vskip+0.1cm

The next lemma may be well-known but we prove it for the completeness.
\vskip+0.1cm

{\bf Lemma E.} {\em Let $(X,\rho, \mu)$ be a non-homogeneous space
with the growth condition. Suppose that $\sigma>-1$. Then there
exists a positive constant $c$ such that for all $a\in X$ and
$r>0$, the inequality

$$ I(a, r, \sigma):= \int_{B(a,r)} \rho(a,x)^{\sigma} d\mu \leq c r^{\sigma+1} $$
holds.}

\vskip+0.1cm

{\bf Proof.} Let $\sigma \geq 0$. Then the result is obvious because of the growth condition for $\mu$.  Further, assume that $-1<\sigma<0$. We have

$$ I(a, r, \sigma) = \int_0^{\infty} \mu \{ x\in B(a, r): \rho(a,x)^{\sigma} >\lambda \} d\lambda $$

$$ = \int_{0}^{\infty} \mu(B(a,r)\cap B(a, \lambda^{1/\sigma})) d\lambda = \int_0^{r^{\sigma}} + \int_{r^{\sigma}}^{\infty} :=
I^{(1)}(a, r, \sigma) + I^{(2)}(a, r, \sigma). $$

By the growth condition for $\mu$ we have

$$ I^{(1)}(a, r, \sigma)\leq r^{\sigma} \mu(B(a,r))\leq c r^{\sigma+1}, $$
while for $I^{(2)}(a, r, \sigma)$ we find that
$$ I^{(2)}(a, r, \sigma) \leq c \int_{r^{\sigma}}^{\infty} \lambda^{1/\sigma} d\lambda = c r^{\sigma+1} $$
because $1/\sigma< -1$.  $\Box$

\vskip+0.1cm

The Following statement is the trace inequality for the operator
$K_{\alpha}$ (see \cite{Adm} for the case of Euclidean spaces  and, e.g.,
\cite{GGKK} or \cite{EKM}, Th. 6.2.1, for an SHT). \vskip+0.1cm

{\bf Theorem F.} {\em Let $(X,\rho, \mu)$ be an SHT. Suppose that
$1<p<q<\infty$ and $0<\alpha<1/p$. Assume that $\nu$ is another
measure on $X$. Then $K_{\alpha}$ is bounded from $L^p(X,\mu)$ to
$L^q(X,\nu)$ if and only if
$$ \nu B  \leq c (\mu B)^{q(1/p-\alpha)}$$
for all balls $B$ in $X$.}

\vskip+0.1cm





\vskip+0.1cm

\section{Main results}

In this section we formulate the main results of the paper.  We begin with the case of an SHT.


\medskip

{\bf Theorem 3.1.} {\em Let $(X, \rho, \mu)$ be an SHT and let
$1<p<q<\infty$. Suppose that  $0<\alpha<1/p$,
$0<\lambda_1<1-\alpha p$ and  $\lambda_2/q= \lambda_1 /p$. Then
$K_{\alpha}$ is bounded from $L^{p,\lambda_1}(X,\mu)$ to $L^{q,\lambda_2}(X,\nu,\mu)$ if and only if there is a positive
constant $c$ such that $$ \nu(B)\leq c \mu(B)^{q(1/p-\alpha)},\eqno{(3)}$$
for all balls $B$.}

\vskip+0.1cm

The next statement is a consequence of Theorem 3.1.
\vskip+0.1cm

{\bf Theorem 3.2.} {\em  Let $(X, \rho, \mu)$ be an SHT and let $1<p<q<\infty$. Suppose that  $0<\alpha <1/p$, $0<\lambda_1<1-\alpha p$ and
$\lambda_2/q = \lambda_1/p$.  Then for the boundedness of $K_\alpha$ from $L^{p,\lambda_1}(X, \mu)$ to $L^{q,\lambda_2}(X, \mu)$ it is necessary and sufficient that $q= p/(1-\alpha p)$. }

\medskip

For non-homogeneous spaces we have the following statements:
\vskip+0.1cm

{\bf Theorem 3.3.} {\em Let $(X,\rho,\mu)$ be a non-homogeneous
space with the growth condition. Suppose that $1< p \le q<\infty$,
$1/p-1/q \le \alpha<1$ and $\alpha\neq 1/p$.  Suppose also that $p
\alpha-1<\beta<p-1$, $0<\lambda_1<\beta-\alpha p +1$ and
$\lambda_1q=\lambda_2p$. Then $\ia$ is bounded from
$M_\beta^{p,\lambda_1}(X,\mu)$ to $M_\gamma^{q,\lambda_2}(X,\mu)$,
where $\gamma=q(1/p+\beta/p-\alpha)-1.$}

\vskip+0.1cm

{\bf Theorem  3.4.} {\em Suppose that $(X,\rho,\mu)$ is a
quasimetric measure space and $\mu$ satisfies condition $(2)$. Let
$1<p<q<\infty$. Assume that $0<\alpha<1$, $0<\lambda_1< p/q$ and
$s\lambda_1/p=\lambda_2/q$.  Then the operator $I_{\alpha}$ is
bounded from $M^{p,\lambda_1 s}(X,\mu)$ to $M^{q,\lambda_2}(X,
\mu)$.}

\vskip+0.1cm


\vskip+0.2cm

\section{Proof of the Main Results}

In this section we give the proofs of the main results.

\medskip

{\em Proof of Theorem} 3.1. {\em Necessity.} Suppose $K_\alpha$ is
bounded from $L^{p,\lambda_1}(\mu)$ to $L^{q,\lambda_2}(X, \nu,
\mu).$ Fix $B_0:= B(x_0,r_0)$. For $x, y\in B_0$, we have that
$$B(x,\rho(x,y))\subseteq B(x, a_1(a_0+1)r_0)\subseteq B(x_0, a_1(1+a_1(a_0+1)) r_0).$$
Hence using the doubling condition for $\mu$, it is easy to see that

$$\mu(B_0)^{\alpha} \leq c K_\alpha \chi_{B_0}(x),\quad x\in B_0.$$

Consequently, using the condition $\lambda_2/q= \lambda_1/p$, the
boundedness of $K_{\alpha}$ from $L^{p,\lambda}(X,\mu)$ to $L^{q,
\lambda_2}(X,\nu, \mu )$ and Lemma D we find that

$$ \mu(B_0)^{\alpha-\lambda_1 /p} \nu(B_0)^{1/q} \leq c \|K_{\alpha}\chi_{B_0}\|_{L^{q,\lambda_2}(X,\nu,\mu)} $$

$$ \leq c \|\chi_{B_0}\|_{L^{p,\lambda_1}(X,\mu)} \leq c \mu(B_0)^{(1-\lambda_1)/p}. $$

Since $c$ does not depend on $B_0$ we have condition $(3)$.

\medskip
{\em Sufficiency.} Let $B:=B(a,r)$, $ \tilde{B}:=B(a,2a_1r)$ and $f\geq 0$.   Write
$f \in L^{p,\lambda_1}(\mu)$ as $f=f_1+f_2:=f\chi_{\tilde
B}+f\chi_{{\tilde B}^{\rm C}},$ where $\chi_B$ is a characteristic
function of $B$.  Then we have
$$ S:= \int_{B} (K_{\alpha}f(x))^q d\nu(x) \leq c \bigg(
\int_{B} (K_{\alpha}f_1(x))^q d\nu(x) + \int_{B}
(K_{\alpha}f_2(x))^q d\nu(x)\bigg) := c(S_1 + S_2). $$

Applying Theorem F and the fact $\mu \in (DC)$ we find that
$$ S_1 \leq  \int_X (K_{\alpha} f_1)^{q}(x) d\nu(x) \leq c
\bigg(\int_{B(a, 2a_1 r)}  (f(x))^p  d\mu(x)\bigg)^{q/p}. $$

Now observe that if $\rho(a,x)<r$ and $\rho(a,y)>2a_1r$, then
$\rho(a,y)>2a_1\rho(a,x)$. Consequently, using the facts $\mu\in
(RDC)$ (see (1)), $0<\lambda_1<1-\alpha p$ and condition (3) we
have

$$ S_2 \leq c  \int_{B(a,r)} \bigg( \int_{\rho(a,y)>r}
\frac{f(y)}{\mu B(a, \rho(a,y)))^{1-\alpha}} d\mu(y) \bigg)^{q}
d\nu(x)  $$

$$  =\nu (B)  \bigg[ \sum_{k=0}^{\infty} \int_{B(a, \eta_1^{k+1} r)\setminus B(a,\eta_1^k r)}
\frac{f(y)}{\mu B(a, \rho(a,y)))^{1-\alpha}}d\mu(y) \bigg]^{q}
$$

$$ \leq c \nu (B) \bigg[ \sum_{k=0}^{\infty} \bigg(\int_{B(a, \eta_1^{k+1}r )} (f(y))^p
d\mu(y)\bigg)^{1/p}$$

$$ \times  \bigg( \int_{B(a, \eta_1^{k+1}
r)\setminus B(a,\eta_1^kr)} \mu B(a, \rho(a,y))^{(\alpha-1)p'}
d\mu(y) \bigg)^{1/p'} \bigg]^{q} $$

$$ \leq c \| f \|^q_{L^{p,\lambda_1}(X,\mu)} \nu (B) \bigg( \sum_{k=0}^{\infty}
\mu B(a, \eta_1^{k+1} r))^{\lambda_1/p+ \alpha-1+1/p'}\bigg)^q
$$

$$ \leq c \| f \|^q_{L^{p,\lambda_1}(X,\mu)} \nu (B)  \mu (B)^{(\lambda_1/p+\alpha-1/p)q} \bigg(\sum_{k=0}^{\infty} \eta_2^{k
(\lambda_1/p+\alpha-1/p)}\bigg)^{q} $$

$$ \leq c \| f\|^q_{L^{p,\lambda_1}(X,\mu)} \mu(B)^{q \lambda_1/p} = c \|f\|^{q}_{L^{p,\lambda_1}(X,\mu)}\mu(B)^{\lambda_2},$$
where the positive constant $c$ does not depend on $B$.  Now the
result follows immediately.  $\Box$

\vskip+0.1cm

{\em Proof of Theorem} 3.2. {\em Sufficency.} Assuming $\alpha= 1/p-1/q$ and $\mu=\nu$ in Theorem 3.1 we have that $K_{\alpha}$ is bounded from
$L^{p,\lambda_1}(X,\mu)$ to $L^{q, \lambda_2}(X, \mu)$.

{\em Necessity.} Suppose that $K_{\alpha}$ is bounded from $L^{p, \lambda_1}(X, \mu)$ to $L^{q, \lambda_2} (X, \mu)$. Then by Theorem 3.1 we have
$$ \mu(B)^{1/q-1/p+\alpha} \leq c. $$

The conditions $\mu(X)= \infty$ and $\mu\{ x\}= 0$ for all $x\in X$ implies that $\alpha= 1/p-1/q$. $\Box$
\vskip+0.1cm

{\em Proof of Theorem} 3.3. Let $f\geq 0$. For $x,\,a\in X,$ let
us introduce the following notation:

$$   E_1(x)\, :=\,
\biggl\{y:\frac{\rho(a,y)}{\rho(a,x)}<\frac{1}{2a_1}\biggr\};\;\;
E_2(x)\, :=\, \biggl\{y:\frac{1}{2a_1}\leq
\frac{\rho(a,y)}{\rho(a,x)}\leq 2a_1\biggr\}; $$
$$ E_3(x)\, :=\,
\biggl\{y:2a_1<\frac{\rho(a,y)}{\rho(a,x)}\biggr\}.$$

For $i=1, 2, 3$,  $r>0$ and $a\in X,$ we denote
$$
S_i:=\int_{\rho(a,x)<r}\rho(a,x)^\gamma\left(\int_{E_i(x)}
f(y)\rho(x,y)^{\alpha-1}\,d\mu(y)\right)^q d\mu(x).
$$

If $y\in E_1(x),$ then $\rho(a,x)<2a_1 a_0 \rho(x,y).$ Hence, it
is easy to see that
$$
S_1\le C\,\int_B\rho(a,x)^{\gamma+
q(\alpha-1)}\left(\int_{\rho(a,y)<\rho(a,x)}f(y)\,d\mu(y)\right)^qd\mu(x).
$$

Taking into account the condition $\gamma<(1-\alpha)q-1$ we have

$$ \int_{ \rho(a,x)>t } \rho(a,x)^{\gamma+ q(\alpha-1)}
d\mu(x) = \sum_{n=0}^{\infty} \int_{B(a, 2^{k+1}t)\setminus
B(a,2^{k}t)} (\rho(a,x))^{\gamma+(\alpha-1)q}d\mu(x) $$

$$\leq c \sum_{n=0}^{\infty} (2^k t)^{\gamma +q(\alpha-1)+1} = c t^{\gamma +q(\alpha-1)+1}, $$
while    the condition $\beta<p-1$ implies
$$ \int_{\rho(a,x)<t}  \rho(a,x)^{\beta( 1-p')+ 1}  d\mu(x) \leq c
t^{\beta(1-p')+1}. $$ Hence

$$ \sup_{a\in X, t>0} \bigg(\int_{\rho(a,x)>t} \rho(a,x)^{\gamma+ q(\alpha-1)}
d\mu(x)\bigg)^{1/q} \bigg(\int_{B(a,t)} \rho(a,y)^{\beta(1-p')}
d\mu(y)\bigg)^{1/p'} <\infty. $$

Now using Theorem C we have

$$ S_1 \le c \, \bigg( \int_B \rho(a,x)^{\beta} (f(y)) \, d\mu(y) \bigg)^{q/p}
\leq c \| f \|^q_{ M^{p,\lambda_1}_{\beta}(X,\mu)} r^{ \lambda_1
q/p}= c \| f \|^q_{ M^{p,\lambda_1}_{\beta}(X,\mu)}
r^{\lambda_2}.$$

Further, observe that if $\rho(a,y)>2 a_1 \rho(a,x)$, then $ \rho
(a,y)\leq a_1 \rho(a,x) +a_1 \rho(a,y) \leq \rho(a,y)/2 +a_1
\rho(x,y). $ Hence $ \rho(a,y)/(2a_1)\leq \rho(x,y)$.
Consequently, using the growth condition for $\mu$ , the fact
$\lambda_1 <\beta-\alpha p+1$ and Lemma E we find that

$$ S_3 \leq c  \int_{B(a,r)} \rho(a,x)^{\gamma} \bigg( \int_{\rho(a,y)> \rho(a,x)}
\frac{f(y)}{\rho(a,y)^{1-\alpha}} d\mu(y) \bigg)^{q} d\mu(x) $$

$$ \leq c  \int_{B(a,r)} \rho(a,x)^{\gamma}  \bigg( \sum_{k=0}^{\infty} \int_{B (a, 2^{k+1} \rho(a,x))\setminus
B(a, 2^k \rho(a,x))} \frac{f(y)}{\rho(a,y)^{1-\alpha}} d\mu(y)
\bigg)^{q} d\mu(x) $$

$$ \leq c  \int_{B(a,r)} \rho(a,x)^{\gamma}  \bigg[ \sum_{k=0}^{\infty} \bigg(
\int_{B(a, 2^{k+1}\rho(a,x))} f^p(y) \rho(a,y)^{\beta}
d\mu(y)\bigg)^{1/p} $$

$$\times \bigg(\int_{B (a, 2^{k+1}
\rho(a,x))\setminus B(a, 2^k \rho(a,x))} \rho(a,y)^{\beta(1-p')+
(\alpha-1)p'} d\mu(y) \bigg)^{1/p'}\bigg]^q  d\mu(x) $$

$$ \leq c \| f \|^q_{ M^{p,\lambda_1}_{\beta}(X,\mu)}
\int_{B(a,r)} \rho(a,x)^{\gamma}  $$

$$ \times \bigg(
\sum_{k=0}^{\infty} (2^k \rho(a,x))^{\lambda_1/p+\alpha-1-\beta/p}
(\mu B(a, 2^{k+1} \rho(a,x)))^{1/p'}\bigg)^{q} d\mu (x) $$

$$ \leq c \| f \|^q_{M^{p,\lambda_1}_{\beta}(X,\mu)}
\int_{B(a,r)} \rho(a,x)^{\gamma}  \bigg( \sum_{k=0}^{\infty} (2^k
\rho(a,x))^{\lambda_1/p+\alpha-1/p-\beta/p} \bigg)^{q} d\mu (x)
$$

$$ \leq c \| f \|^q_{ M^{p,\lambda_1}_{\beta}(X,\mu)}
\int_{B(a,r)} \rho(a,x)^{(\lambda_1/p+ \alpha-1/p-\beta/p)q
+\gamma} d\mu(x) $$

$$ = c  \| f \|^q_{ M^{p,\lambda_1}_{\beta}(X,\mu)}
\int_{B(a,r)} \rho(a,x)^{\lambda_1 q/p -1} d\mu(x)  \leq c \| f
\|^q_{ M^{p,\lambda_1}_{\beta}(X,\mu)} r^{ \lambda_1 q/p}   $$
$$ =  c \| f \|^q_{ M^{p,\lambda_1}_{\beta}(X,\mu)}
r^{ \lambda_2}.  $$

So, we conclude that
$$ S_3 \leq c \| f \|^q_{ M^{p,\lambda_1}_{\beta}(X,\mu)}
r^{ \lambda_2}. $$

To estimate $S_2$ we consider two cases. First assume that
$\alpha<1/p$. Let

$$ E_{k,r}:= \{ x: 2^k r \leq  \rho(a,x) < 2^{k+1} r \}; $$

$$ F_{k,r}:= \{ x: 2^{k-1} r/ a_1  \leq \rho(a,x)< a_1 2^{k+2} r \}. $$

Assume that $p^{*}=p/(1-\alpha p)$. By H\"older's
inequality,  Corollary B and the assumption $\gamma=  q(1/p+\beta/p -\alpha)-1$ we have

$$ S_2 = \sum_{k=-\infty}^{-1} \int_{E_{k,r}}\rho(a,x)^{\gamma}\bigg( \int_{E_2(x)} f(y)
\rho(x,y)^{\alpha-1} d\mu(y) \bigg)^q d\mu(x)  $$

$$ \leq \sum_{k=-\infty}^{-1}
\bigg(\int_{E_{k,r}}\rho(a,x)^{\gamma}\bigg( \int_{E_2(x)} f(y)
\rho(x,y)^{\alpha-1} d\mu(y) \bigg)^{p^*} d\mu(x)\bigg)^{q/p'}$$
$$\times \bigg( \int_{E_{k,r}} \rho(a,x)^{\gamma p^*/(p^* -q)}
d\mu(x)\bigg)^{(p^*-q)/p^*} $$

$$ \leq c \sum_{k=-\infty}^{-1} 2^{k(\gamma+(p^*-q)/p^*)}
\bigg(\int_{X} I_{\alpha} ( f\chi_{F_{k,r}})(x))^{p^*} d\mu(x)
\bigg)^{q/p^*} $$

$$ \leq  c \sum_{k=-\infty}^{-1}
2^{k(\gamma+(p^*-q)/p^*)} \bigg( \int_{F_{k,r}} (f(x))^{p} d\mu(x)
\bigg)^{q/p} \leq c \bigg(\int_{B(a, 2a_1 r)} \rho(a,x)^{\beta} (f(x))^{p}
d\mu(x) \bigg)^{q/p} $$

$$ \leq c \| f \|^q_{ M^{p,\lambda_1}_{\beta}(X,\mu)}
r^{ \lambda_1 q/p} =  c \| f \|^q_{ M^{p,\lambda_1}_{\beta}(X,\mu)}
r^{ \lambda_2 }. $$

Let us now consider the case $1/p<\alpha<1$.

First notice that (see \cite{KM1})

$$ \int_{E_2(x)}(\rho(x,y)^{(\alpha-1)p'} d\mu(y)  \leq c
\rho(a,x)^{1+(\alpha-1)p'}, $$ where the positive constant $c$
does not depend on $a$ and $x$.

This estimate and  H\"older's inequality yield

$$ S_2 \leq c \sum_{k=-\infty}^{-1}
\bigg(\int_{E_{k,r}} \rho(a,x)^{\gamma+[(\alpha-1)p'+1)]q/p'}
\bigg( \int_{E_2(x)} (f(y))^p d\mu(y) \bigg)^{q/p}
d\mu(x)\bigg)^{q/p'}
$$

$$ \leq c \sum_{k=-\infty}^{-1}
\bigg(
\int_{E_{k,r}}\rho(a,x)^{\gamma+[(\alpha-1)p'+1)]q/p'}d\mu(x)\bigg)
\bigg(\int_{F_{k,r}} (f(y))^{p} d\mu(y) \bigg)^{q/p} $$

$$ \leq c \sum_{k=-\infty}^{-1} (2^k r)^{\gamma+[(\alpha-1)p'+1)]q/p'+1}
\bigg(\int_{F_{k,r}} (f(y))^{p} d\mu(y) \bigg)^{q/p} $$

$$ = c \sum_{k=-\infty}^{-1} 2^{k \beta q/ p}
\bigg(\int_{F_{k,r}} (f(y))^{p} d\mu(y) \bigg)^{q/p} \leq c \bigg(
\int_{B(a, 2a_1 r)} (f(y))^{p} \rho(a,y)^{\beta} d\mu(y) \bigg)^{q/p} $$

$$ \leq c \| f \|^q_{ M^{p,\lambda_1}_{\beta}(X,\mu)}
r^{ \lambda_1 q/p} =  c \| f \|^q_{ M^{p,\lambda_1}_{\beta}(X,\mu)}
r^{ \lambda_2 }. $$
Now the result follows immediately. $\Box$

\vskip+0.1cm

{\em Proof of Theorem} 3.4. Let $f\geq 0$. Suppose that $a\in X$ and
$r>0$. Suppose also that $f_1= f\chi_{B(a, 2a_1 r)}$ and $f_2=f-f_1$. Then  $I_{\alpha}f = I_{\alpha}f_1 + I_{\alpha}f_2$.  Consequently,
$$ \int_{B(a,r)}(I_{\alpha}f(x))^q d\mu(x) \leq 2^{q-1} \bigg( \int_{B(a,r)}(I_{\alpha}f_1(x))^q d\mu(x) $$
$$ + \int_{B(a,r)}(I_{\alpha}f_2(x))^q d\mu(x) \bigg) := 2^{q-1}( S^{(1)}_{a,r}+ S^{(2)}_{a,r}). $$
Due to Theorem A and the condition $s\lambda_1/p= \lambda_2/q$ we have

$$ S^{(1)}_{a,r} \leq c \bigg(\int_{B(a, 2a_1 r)} (f(x))^p d\mu(x)\bigg)^{q/p} $$
$$= c \bigg( \frac{1}{(2a_1 r)^{\lambda_1 s}} \int_{B(a, 2a_1 r)} (f(x))^p dx \bigg)^{q/p} r^{\lambda_1 s q/p} \leq c \| f\|^q_{M^{p,\lambda_1s}(X,\mu)} r^{\lambda_2}. $$

Now observe that if $x\in B(a, r)$ and $y\in X \setminus B(a, 2a_1
r)$, then $ \frac{\rho(a, y)}{2a_1} \leq \rho(x,y). $ Hence
H\"older's inequality, condition (2) and the condition
$0<\lambda_1<p/q$ yield

$$ I_{\alpha}f_2(x)  = \int_{X\setminus B(a,2a_1r)} f(y)/\rho(x,y)^{1-\alpha} d\mu(y) $$

$$= \sum_{k=0}^{\infty} \bigg( \int_{B(a, 2^{k+2}a_1 r)\setminus B(a, 2^{k+1}a_1 r)} (f(y))^p
d\mu(y)\bigg)^{1/p} $$

$$\times \bigg( \int_{B(a, 2^{k+2}a_1 r)\setminus B(a, 2^{k+1}a_1
r)} \rho(a,y)^{(\alpha-1)p'}d\mu(y)\bigg)^{1/p'} $$

$$ \leq c \sum_{k=0}^{\infty} \bigg( \frac{1}{(2^{k+1}a_1 r)^{\lambda_1 s}} \int_{B(a, 2^{k+1}a_1 r)} (f(y))^p d\mu(y) \bigg)^{1/p}(2^{k}a_1r)^{\lambda_1 s/p +\alpha -1 + s/p'} $$

$$ \leq c \| f\|_{M^{p,\lambda_1 s}(X,\mu)}
r^{\lambda_1s /p+\alpha-1+s/p'}.$$

Consequently, by the assumptions $s\lambda_1/p =\lambda_2/q$ and $s=\frac{pq(1-\alpha)}{pq+p-q}$ we conclude that

$$ S^{(2)}_{a,r}\leq c \| f\|^q_{M^{p,\lambda_1 s}(X,\mu)}
r^{(\lambda_1 s/p+\alpha-1+s/p')q +s}= c \|
f\|^q_{M^{p,\lambda_1 s}(X,\mu)} r^{\lambda_2}. $$

Summarazing the estimates derived above we finally have the desired result. $\Box$

\vskip+0.5cm





\vskip+0.5cm
\begin{center}
{\bf Acknowledgement}
\end{center}
\vskip+0.4cm

The second and third authors were partially
supported by the INTAS Grant  No. 05-1000008-8157 and the Georgian National Science Foundation Grant No. GNSF/ST07/3-169.
\vskip+1cm

\end{document}